# Construction of fuzzy valued recurrent fractal interpolation functions and their properties

Hyang Choe, MiGyong Ri*, CholHui Yun

Faculty of Mathematics, **Kim Il Sung** University, Pyongyang, Democratic People's Republic of Korea

*mk.ri1004@ryongnamsan.edu.kp

ch.yun@ryongnamsan.edu.kp

Abstract. In the process of measuring objects with local self-similarity, such as satellite images or coastlines, we obtain a data set with both local self-similarity and uncertainty. To better interpolate such data sets, an interpolation function with both local self-similarity and uncertainty is necessary.

In this paper, we propose a construction of fuzzy valued recurrent fractal interpolation function using recurrent iterated function system that interpolates the given data set of fuzzy numbers. And we show some properties of the constructed fuzzy valued RFIFs: Hölder continuity and stability of the interpolation function due to perturbations in the data set or the vertical scaling factors.

## 1. Introduction

In nature, there are many objects with very complex structures, which are difficult to mod el with classical interpolation functions, such as coastline, seismic, image, meteorological data and so on. To better model these objects mathematically, Barnsley([3]) introduced a fractal inte rpolation function(FIF) by using an iterated function system(IFS). FIF models objects with self -similarity better than classical interpolation functions. Later, recurrent fractal interpolation funct ions(RFIF) using recurrent iterated function systems(RIFS) were introduced to better model obj ects with local self-similarity, such as image data or coastline([4,7,20]). Also, construction and properties of various kinds of FIFs, such as hidden variable FIF and hidden variable recurrent FIF, have been studied([5,6,9,11,14,18,22]).

The real data, given by experiments, observations, measurements, etc., are perturbed by various factors and so, they have uncertainty. Hence, in data processing and analysis, it is necessary to model the data with uncertainty accurately. To express the degree of uncertainty mathematically, Zadeh([23]) introduced the concept of fuzzy sets. To interpolate the data set of fuzzy numbers, fuzzy interpolation polynomial of Lagrange type by Lowen([15]), fuzzy Lagrange polynomial and fuzzy spline by Kaleva([13]), and fuzzy Hermite interpolation polynomial by Goghary([12]) had been studied. Abbasbandy et al.([1]) studied interpolation of fuzzy data by a new set of spline functions.

The real data observed from highly irregular phenomena have uncertainty due to various factors, so these data have not only fractal characteristic but also uncertainty. Therefore, interpolation representing both local self-similarity and uncertainty is necessary. That's why, many mathematicians have attempted to apply the concept of fuzzification introduced by Zadeh to fractals([2,8,16,17,19]). Xiao et al.([19]) proposed a method to model data with uncertainty using FIF. They converted the points of the given data into fuzzy numbers and then generated a family of FIFs using sets of starting points and endpoints of level sets of fuzzy numbers. Yun et al. ([21]) generalized iterated fuzzy set system to one with recurrent structure and constructed a recurrent fuzzy fractal. However, unfortunately, the obtained interpolation functions are not fuzzy valued interpolation functions with self-similarity.

In this paper, we will construct fuzzy valued RFIFs by using RIFSs to better describe natural objects with both local self-similarity and uncertainty and their properties.

The remainder of this paper is organized as follows: In Section 2, we propose construction of a fuzzy valued RFIF interpolating a given data set of fuzzy numbers by using RIFS. In Section 3, we study properties of the constructed fuzzy valued RFIFs: Hölder continuity and stability in perturbations of the data set or the vertical scaling factors.

## 2. Construction of recurrent iterated function system and fuzzy valued recurrent fractal interpolaion function from fuzzy value data

Let's $\mathbf{R_F}$ be the space of fuzzy numbers over $\mathbf{R}$. For a fuzzy number $u \in \mathbf{R_F}$ and $0 \le \lambda \le 1$, $[u]^\lambda$, the level set and its length are denoted by $[u^-(\lambda),\ u^+(\lambda)]$ and $len([u]^\lambda)$, respectively.

$(\mathbf{R_F}, d_\infty)$ is a complete metric space with respect to the following Hausdorff metric $d_\infty$ on $\mathbf{R_F}$ ([10]):

$$d_\infty(u,\ v) = \sup_{\lambda\in[0,1]} \max\{|u^-(\lambda) - v^-(\lambda)|,\ |u^+(\lambda) - v^+(\lambda)|\},\ u,\ v \in \mathbf{R_F}.$$

$C(\mathrm{K},\ \mathbf{R_F})$, the space of continuous fuzzy valued functions on a compact subset K of $\mathbf{R}$ is a complete metric space with respect to the following metric $D$ ([10]):

$$D(f,\ g) = \sup_{x\in K} d_\infty(f(x),\ g(x)),\ f,\ g \in C(\mathrm{K},\ \mathbf{R_F}).$$

Moreover, for any $u_i \in \mathbf{R_F},\ i = 0,\ 1,\ \cdots,\ n$, the following $C^*(I,\ \mathbf{R_F})$ is also a complete metric space with respect to the metric $D$:

$$C^*(I,\ \mathbf{R_F}) = \{\varphi : I \to \mathbf{R_F} \mid \varphi \text{ is continuous and } \varphi(x_i) = u_i,\ i = 0,\ 1,\ \cdots,\ n\}$$

Through this paper, we denote the Lipschitz coefficient of a Lipschitz function $f$ by $L_f$ and contractive coefficient of a contractive function $f$ by $c_f$, Hölder coefficient of a Hölder continuous function $f$ by $H_f$. For a fuzzy valued function $\mathbf{f} \in C(\mathrm{K}, \mathbf{R_F})$ and any $\lambda \in [0, 1]$, $\{(x,\ (\mathbf{f}(x))^-(\lambda)) \mid x \in \mathrm{K}\} \cup \{(x,\ (\mathbf{f}(x))^+(\lambda)) \mid x \in \mathrm{K}\}$ is called the $\lambda$-level set of $\mathbf{f}$.

### 1) Construction of recurrent iterated function system

A data set of fuzzy numbers is given as follows:

$$\mathbf{P} = \{(x_i,\ u_i) \in \mathbf{R} \times \mathbf{R_F} \mid i = 0,\ 1, \cdots,\ n\}\ (x_0 < x_1 < \cdots < x_n).$$

Let's $I = [x_0,\ x_n]$ and $I_i = [x_{i-1},\ x_i],\ i = 1,\ \cdots,\ n$. For every $I_i$, we correspond an interval

$\widetilde{I}_{\sigma(i)} = [x_{s_i},\ x_{e_i}]\ (e_i - s_i \ge 2)$ consisting of several subintervals to $I_i$.

The functions $\widetilde{l}_i : \widetilde{I}_{\sigma(i)} \to I_i,\ i = 1,\ \cdots,\ n$ are defined as follows:

$$\widetilde{l}_i(x) = \frac{x_i - x_{i-1}}{x_{e_i} - x_{s_i}}(x - x_{s_i}) + x_{i-1},\ x \in \widetilde{I}_{\sigma(i)}, \tag{1}$$

where $c_{\widetilde{l}_i} = \dfrac{x_i - x_{i-1}}{x_{e_i} - x_{s_i}} < 1$.

The mappings $\widetilde{\mathbf{F}}_i : \widetilde{I}_{\sigma(i)} \times \mathbf{R_F} \to \mathbf{R_F},\ i = 1,\ \cdots,\ n$ are defined as follows:

$$\widetilde{\mathbf{F}}_i(x,\ u) = \alpha_i u \oplus \widetilde{\mathbf{q}}_i(\widetilde{l}_i(x)),\ (x,\ u) \in \widetilde{I}_{\sigma(i)} \times \mathbf{R_F}, \tag{2}$$

where $0 \le \alpha_i < 1$ and $\widetilde{\mathbf{q}}_i : I_i \to \mathbf{R_F},\ i = 1,\ \cdots,\ n$ are Lipschitz mappings satisfying the following conditions:

$$\widetilde{\mathbf{F}}_i(x_{s_i},\ u_{s_i}) = u_{i-1},\ \widetilde{\mathbf{F}}_i(x_{e_i},\ u_{e_i}) = u_i. \tag{3}$$

Now we give an example of Lipschitz mappings $\widetilde{\mathbf{q}}_i$ satisfying (3).

**Example 1.** The data set $\mathbf{P}$ is given and we choose the vertical scaling factors $\alpha_i,\ i = 1,\ \cdots,\ n$ such that they are satisfied the following conditions:

(a1) $\forall \lambda \in [0,\ 1],\ len([u_{i-1}]^\lambda) \ge len([\alpha_i u_{s_i}]^\lambda),\ len([u_i]^\lambda) \ge len([\alpha_i u_{e_i}]^\lambda)$.

(a2) $u_{i-1}^-(\lambda) - \alpha_i u_{s_i}^-(\lambda)$ and $u_i^-(\lambda) - \alpha_i u_{e_i}^-(\lambda)$ are increasing with respect to $\lambda$.

(a3) $u_{i-1}^+(\lambda) - \alpha_i u_{s_i}^+(\lambda)$ and $u_i^+(\lambda) - \alpha_i u_{e_i}^+(\lambda)$ are decreasing with respect to $\lambda$.

The following mappings $\widetilde{\mathbf{q}}_i,\ i = 1,\ \cdots,\ n$ are Lipschitz ones satisfying (3):

$$\widetilde{\mathbf{q}}_i(x) = \left( \frac{x - x_{i-1}}{x_i - x_{i-1}} u_i \oplus \frac{x - x_i}{x_{i-1} - x_i} u_{i-1} \right) \vee_H \left[ \alpha_i \left( \frac{\widetilde{l}_i^{-1}(x) - x_{s_i}}{x_{e_i} - x_{s_i}} u_{e_i} \oplus \frac{\widetilde{l}_i^{-1}(x) - x_{e_i}}{x_{s_i} - x_{e_i}} u_{s_i} \right) \right],\ x \in I_i, \tag{4}$$

where $\vee_H$ is Hukuhara Difference and the conditions (a1)-(a3) guarantee the existence of Hukuhara Difference in (4).

Firstly, $\widetilde{\mathbf{q}}_i,\ i = 1,\ \cdots,\ n$ satisfy (3) because

$$\widetilde{\mathbf{F}}_i(x_{s_i},\ u_{s_i})=\alpha_i u_{s_i}\oplus\widetilde{\mathbf{q}}_i(\widetilde{l}_i(x_{s_i}))=\alpha_i u_{s_i}\oplus\widetilde{\mathbf{q}}_i(x_{i-1})=\alpha_i u_{s_i}\oplus[u_{i-1}\vee_H(\alpha_i u_{s_i})]=u_{i-1},$$

$$\widetilde{\mathbf{F}}_i(x_{e_i},\ u_{e_i})=\alpha_i u_{e_i}\oplus\widetilde{\mathbf{q}}_i(\widetilde{l}_i(x_{e_i}))=\alpha_i u_{e_i}\oplus\widetilde{\mathbf{q}}_i(x_i)=\alpha_i u_{e_i}\oplus[u_i\vee_H(\alpha_i u_{e_i})]=u_i.$$

Secondly, $\widetilde{\mathbf{q}}_i: I_i\to\mathbf{R}_\mathbf{F}$ are Lipschitz mappings. For any $x,\ y\in I_i$ and $0\le\lambda\le 1$, we have

$$|(\widetilde{\mathbf{q}}_i(x))^-(\lambda)-(\widetilde{\mathbf{q}}_i(y))^-(\lambda)|=$$

$$=\left|\left\{\left(\frac{x-x_{i-1}}{x_i-x_{i-1}}u_i\oplus\frac{x-x_i}{x_{i-1}-x_i}u_{i-1}\right)^-(\lambda)-\left[\alpha_i\left(\frac{\widetilde{l}_i^{-1}(x)-x_{s_i}}{x_{e_i}-x_{s_i}}u_{e_i}\oplus\frac{\widetilde{l}_i^{-1}(x)-x_{e_i}}{x_{s_i}-x_{e_i}}u_{s_i}\right)\right]^-(\lambda)\right\}\right.$$

$$\left.-\left\{\left(\frac{y-x_{i-1}}{x_i-x_{i-1}}u_i\oplus\frac{y-x_i}{x_{i-1}-x_i}u_{i-1}\right)^-(\lambda)-\left[\alpha_i\left(\frac{\widetilde{l}_i^{-1}(y)-x_{s_i}}{x_{e_i}-x_{s_i}}u_{e_i}\oplus\frac{\widetilde{l}_i^{-1}(y)-x_{e_i}}{x_{s_i}-x_{e_i}}u_{s_i}\right)\right]^-(\lambda)\right\}\right|$$

$$=\left|\left[\frac{x-x_{i-1}}{x_i-x_{i-1}}(u_i^-(\lambda)-\alpha_i u_{e_i}^-(\lambda))+\frac{x-x_i}{x_{i-1}-x_i}(u_{i-1}^-(\lambda)-\alpha_i u_{s_i}^-(\lambda))\right]\right.$$

$$\left.-\left[\frac{y-x_{i-1}}{x_i-x_{i-1}}(u_i^-(\lambda)-\alpha_i u_{e_i}^-(\lambda))+\frac{y-x_i}{x_{i-1}-x_i}(u_{i-1}^-(\lambda)-\alpha_i u_{s_i}^-(\lambda))\right]\right|$$

$$=\left|\frac{x-y}{x_i-x_{i-1}}(u_i^-(\lambda)-\alpha_i u_{e_i}^-(\lambda))+\frac{x-y}{x_{i-1}-x_i}(u_{i-1}^-(\lambda)-\alpha_i u_{s_i}^-(\lambda))\right|$$

$$\le\frac{|x-y|}{x_i-x_{i-1}}\cdot|[u_i^-(\lambda)-\alpha_i u_{e_i}^-(\lambda)]-[u_{i-1}^-(\lambda)-\alpha_i u_{s_i}^-(\lambda)]|.$$

Since by the condition (a2),

$$u_{i-1}^-(0)-\alpha_i u_{s_i}^-(0)\le u_{i-1}^-(\lambda)-\alpha_i u_{s_i}^-(\lambda)\le u_{i-1}^-(1)-\alpha_i u_{s_i}^-(1),$$

$$u_i^-(0)-\alpha_i u_{e_i}^-(0)\le u_i^-(\lambda)-\alpha_i u_{e_i}^-(\lambda)\le u_i^-(1)-\alpha_i u_{e_i}^-(1),$$

we get

$$|(\widetilde{\mathbf{q}}_i(x))^-(\lambda)-(\widetilde{\mathbf{q}}_i(y))^-(\lambda)|\le$$

$$\le\frac{|x-y|}{x_i-x_{i-1}}\max\{|[u_i^-(1)-\alpha_i u_{e_i}^-(1)]-[u_{i-1}^-(0)-\alpha_i u_{s_i}^-(0)]|,\ |[u_i^-(0)-\alpha_i u_{e_i}^-(0)]-[u_{i-1}^-(1)-\alpha_i u_{s_i}^-(1)]|\}.$$

Similarly, we have

$$|(\widetilde{\mathbf{q}}_i(x))^+(\lambda)-(\widetilde{\mathbf{q}}_i(y))^+(\lambda)|\le$$

$$\le\frac{|x-y|}{x_i-x_{i-1}}\max\{|[u_i^+(1)-\alpha_i u_{e_i}^+(1)]-[u_{i-1}^+(0)-\alpha_i u_{s_i}^+(0)]|,\ |[u_i^+(0)-\alpha_i u_{e_i}^+(0)]-[u_{i-1}^+(1)-\alpha_i u_{s_i}^+(1)]|\}.$$

Hence, we obtain

$$d_\infty(\widetilde{\mathbf{q}}_i(x),\ \widetilde{\mathbf{q}}_i(y)) = \sup_{0\le\lambda\le1} \max\{|(\widetilde{\mathbf{q}}_i(x))^-(\lambda) - (\widetilde{\mathbf{q}}_i(y))^-(\lambda)|,\ |(\widetilde{\mathbf{q}}_i(x))^+(\lambda) - (\widetilde{\mathbf{q}}_i(y))^+(\lambda)|\}$$

$$\le L_{\widetilde{\mathbf{q}}_i} |x-y|,$$

where

$$L_{\widetilde{\mathbf{q}}_i} = \frac{1}{x_i - x_{i-1}} \max\{|[u_i^+(1) - \alpha_i u_{e_i}^+(1)] - [u_{i-1}^+(0) - \alpha_i u_{s_i}^+(0)]|,\ |[u_i^+(0) - \alpha_i u_{e_i}^+(0)] - [u_{i-1}^+(1) - \alpha_i u_{s_i}^+(1)]|,$$
$$|[u_i^-(1) - \alpha_i u_{e_i}^-(1)] - [u_{i-1}^-(0) - \alpha_i u_{s_i}^-(0)]|,\ |[u_i^-(0) - \alpha_i u_{e_i}^-(0)] - [u_{i-1}^-(1) - \alpha_i u_{s_i}^-(1)]|\}. \tag{5}$$

That is, $\widetilde{\mathbf{q}}_i,\ i=1,\ \cdots,\ n$ are Lipschitz mappings satisfying (3).

For any $\theta > 0$, let's define a metric $d_\theta : (I \times \mathbf{R_F}) \times (I \times \mathbf{R_F}) \to \mathbf{R}$ as follows:

$$d_\theta((x,\ u),\ (y,\ v)) = |x-y| + \theta d_\infty(u,\ v),\ (x,\ u),\ (y,\ v) \in I \times \mathbf{R_F}.$$

Then, $(I \times \mathbf{R_F}, d_\theta)$ is a complete metric space. In fact, as we have

$$d_\theta((x_n,\ u_n),\ (x_m,\ u_m)) = |x_n - x_m| + \theta d_\infty(u_n,\ u_m) \ge |x_n - x_m|,$$

$$d_\theta((x_n,\ u_n),\ (x_m,\ u_m)) \ge \theta d_\infty(u_n,\ u_m)$$

for any Cauchy sequence $\{(x_n,\ u_n)\} \subset I \times \mathbf{R_F}$, $\{x_n\}$ and $\{u_n\}$ are Cauchy sequences in the complete metric space $I$ and $\mathbf{R_F}$, respectively. Hence, $\{x_n\}$ and $\{u_n\}$ converge to $x_0$ and $u_0$, respectively. Since

$$d_\theta((x_n,\ u_n),\ (x_0,\ u_0)) = |x_n - x_0| + \theta d_\infty(u_n,\ u_0)$$
$$\le (1+\theta) \cdot \max\{|x_n - x_0|,\ d_\infty(u_n,\ u_0)\},$$

the sequence $\{(x_n,\ u_n)\}$ converges to $(x_0,\ u_0) \in I \times \mathbf{R_F}$.

Let's define the mappings $\widetilde{\mathbf{w}}_i : \widetilde{I}_{\sigma(i)} \times \mathbf{R_F} \to I_i \times \mathbf{R_F},\ i=1,\ \cdots,\ n$ as follows:

$$\widetilde{\mathbf{w}}_i(x,\ u) = (\widetilde{l}_i(x),\ \widetilde{\mathbf{F}}_i(x,\ u)),\ (x,\ u) \in \widetilde{I}_{\sigma(i)} \times \mathbf{R_F}.$$

If $0 < \theta < \dfrac{1 - \max_{i=1,\cdots,n}\{c_{\widetilde{l}_i}\}}{\max_{i=1,\cdots,n}\{L_{\widetilde{\mathbf{q}}_i} \cdot c_{\widetilde{l}_i}\}}$, then the following result holds:

***Theorem 1.*** $\tilde{\mathbf{w}}_i$, $i = 1, \cdots, n$ ***are contractive mappings with respect to the metric*** $d_\theta$ ***.***

Proof. For any $(x, u)$ and $(y, v)$ in $\tilde{I}_{\sigma(i)} \times \mathbf{R_F}$, we get

$$
\begin{aligned}
d_\theta(\tilde{\mathbf{w}}_i(x, u), \tilde{\mathbf{w}}_i(y, v)) &= d_\theta((\tilde{l}_i(x), \tilde{\mathbf{F}}_i(x, u)), (\tilde{l}_i(y), \tilde{\mathbf{F}}_i(y, v)) \\
&= |\tilde{l}_i(x) - \tilde{l}_i(y)| + \theta d_\infty(\tilde{\mathbf{F}}_i(x, u), \tilde{\mathbf{F}}_i(y, v)) \\
&\le c_{\tilde{l}_i} |x - y| + \theta d_\infty(\alpha_i u \oplus \tilde{\mathbf{q}}_i(\tilde{l}_i(x)), \alpha_i v \oplus \tilde{\mathbf{q}}_i(\tilde{l}_i(y))) \\
&\le c_{\tilde{l}_i} |x - y| + \theta[d_\infty(\alpha_i u, \alpha_i v) + d_\infty(\tilde{\mathbf{q}}_i(\tilde{l}_i(x)), \tilde{\mathbf{q}}_i(\tilde{l}_i(y)))] \\
&\le c_{\tilde{l}_i} |x - y| + \theta \cdot \alpha_i d_\infty(u, v) + \theta L_{\tilde{\mathbf{q}}_i} c_{\tilde{l}_i} |x - y| \\
&= (c_{\tilde{l}_i} + \theta L_{\tilde{\mathbf{q}}_i} c_{\tilde{l}_i}) |x - y| + \alpha_i \theta d_\infty(u, v) \\
&\le c_{\tilde{\mathbf{w}}_i}(|x - y| + \theta d_\infty(u, v)) = c_{\tilde{\mathbf{w}}_i} d_\theta((x, u), (y, v)),
\end{aligned}
$$

where $c_{\tilde{\mathbf{w}}_i} = \max\{c_{\tilde{l}_i} + \theta L_{\tilde{\mathbf{q}}_i} c_{\tilde{l}_i}, \alpha_i\}$. Since

$$
\begin{aligned}
c_{\tilde{l}_i} + \theta L_{\tilde{\mathbf{q}}_i} c_{\tilde{l}_i} &< c_{\tilde{l}_i} + \frac{1 - \max\limits_{i=1,\cdots,n}\{c_{\tilde{l}_i}\}}{\max\limits_{i=1,\cdots,n}\{L_{\tilde{\mathbf{q}}_i} \cdot c_{\tilde{l}_i}\}} L_{\tilde{\mathbf{q}}_i} c_{\tilde{l}_i} = c_{\tilde{l}_i} + \frac{L_{\tilde{\mathbf{q}}_i} c_{\tilde{l}_i}}{\max\limits_{i=1,\cdots,n}\{L_{\tilde{\mathbf{q}}_i} \cdot c_{\tilde{l}_i}\}}(1 - \max\limits_{i=1,\cdots,n}\{c_{\tilde{l}_i}\}) \\
&\le c_{\tilde{l}_i} + 1 - \max\limits_{i=1,\cdots,n}\{c_{\tilde{l}_i}\} \le 1
\end{aligned}
$$

and $\alpha_i < 1$, $c_{\tilde{\mathbf{w}}_i} < 1$, $\tilde{\mathbf{w}}_i$ is a contractive mapping with the contractive coefficient $c_{\tilde{\mathbf{w}}_i}$. □

Let a row stochastic matrix $M = (p_{st})_{n \times n}$ be defined as:

$$
p_{st} = \begin{cases} 1/a_s, & I_s \subseteq \tilde{I}_{\sigma(t)}, \\ 0, & I_s \not\subset \tilde{I}_{\sigma(t)}, \end{cases}
$$

where $a_s$ is the number of $\tilde{I}_{\sigma(i)}$ involving $I_s$. If $M$ is irreducible, the following is a hyperb olic recurrent iterated function system(RIFS):

$$
\{I \times \mathbf{R_F};\ M;\ \tilde{\mathbf{w}}_i,\ i = 1, \cdots, n\}. \tag{6}
$$

**2) Construction of fuzzy valued recurrent fractal interpolation function**

We denote the set of all nonempty compact subsets of $I \times \mathbf{R_F}$ by $H(I \times \mathbf{R_F})$ and the product space $H(I\times\mathbf{R_F})^n = H(I\times\mathbf{R_F})\times\cdots\times H(I\times\mathbf{R_F})$ by $\widetilde{H}(I\times\mathbf{R_F})$. Then $(\widetilde{H}(I\times\mathbf{R_F}),\ \widetilde{h})$ is a complete metric space([7]):

$$\widetilde{h}((A_1,\ \cdots,\ A_n),\ (B_1,\ \cdots,\ B_n)) = \max_{i=1,\cdots,n}\{\sup_{x\in A_i}\inf_{y\in B_i} d_\theta(x,\ y),\ \sup_{y\in B_i}\inf_{x\in A_i} d_\theta(x,\ y)\},$$

$$(A_1,\ \cdots,\ A_N),\ (B_1,\ \cdots,\ B_N)\in\widetilde{H}(I\times\mathbf{R_F}).$$

Let's define the transformation $\widetilde{\mathbf{W}}:\widetilde{H}\to\widetilde{H}$ as follows:

$$\widetilde{\mathbf{W}}(\boldsymbol{B}) = \begin{pmatrix} \bigcup_{j\in\Lambda(1)}\widetilde{\mathbf{w}}_1(B_j) \\ \bigcup_{j\in\Lambda(2)}\widetilde{\mathbf{w}}_2(B_j) \\ \cdots \\ \bigcup_{j\in\Lambda(n)}\widetilde{\mathbf{w}}_n(B_j) \end{pmatrix},\quad \boldsymbol{B}=(B_1,\cdots,B_n)\in\widetilde{H}(I\times\mathbf{R_F}),$$

where $\Lambda(i)=\{j;\ c_{ij}=1\},\ i=1,\ 2,\ \cdots,\ n$. Then, $\widetilde{\mathbf{W}}$ is contractive on $\widetilde{H}(I\times\mathbf{R_F})$ and there exists a unique $\boldsymbol{A}=(A_1,\ \cdots,\ A_n)\in\widetilde{H}(I\times\mathbf{R_F})$ such that $\widetilde{\mathbf{W}}(\boldsymbol{A})=\boldsymbol{A}$, i.e., $A_i=\bigcup_{j\in\Lambda(i)}\widetilde{\mathbf{w}}_i(A_j)$. $A=\bigcup_{i=1}^{n}A_i$ is called the attractor of the RIFS $\{I\times\mathbf{R_F};\ M;\ \widetilde{\mathbf{w}}_i,\ i=1,\ \cdots,\ n\}$ and the continuous fuzzy valued function, whose graph is $A$, is called the fuzzy valued recurrent fractal interpolation function(RFIF) interpolating the data set of fuzzy numbers $\mathbf{P}$.

***Theorem 2.*** ***$A$, the attractor of RIFS (6) is the graph of a continuous fuzzy valued fun ction $\widetilde{\mathbf{f}}$.***

Proof. Let's define the mapping $\widetilde{T}$ on $C^*(I,\ \mathbf{R_F})$ as follows:

$$(\widetilde{T}\varphi)(x)=\widetilde{\mathbf{F}}_i(\widetilde{l}_i^{-1}(x),\ \varphi(\widetilde{l}_i^{-1}(x))),\ \varphi\in C^*(I,\ \mathbf{R_F}),\ x\in I_i.$$

$\widetilde{T}:C^*(I,\ \mathbf{R_F})\to C^*(I,\ \mathbf{R_F})$ is well defined. In fact, by the definition of $\widetilde{\mathbf{F}}_i$, $\widetilde{T}\varphi$ is continuous in $I_i$ and we have

$$(\widetilde{T}\varphi)(x_i)=\widetilde{\mathbf{F}}_i(\widetilde{l}_i^{-1}(x_i),\ \varphi(\widetilde{l}_i^{-1}(x_i)))=\widetilde{\mathbf{F}}_i(x_{e_i},\ \varphi(x_{e_i}))=\widetilde{\mathbf{F}}_i(x_{e_i},\ u_{e_i})=u_i,$$

$$(\widetilde{T}\varphi)(x_{i-1})=u_{i-1}.$$

We show that $\widetilde{T}$ is contractive. In fact, for any $\varphi_1$ and $\varphi_2$ in $C^*(I,\ \mathbf{R_F})$, the following holds:

$$D(\widetilde{T}\varphi_1,\ \widetilde{T}\varphi_2) = \sup_{x\in I} d_\infty((\widetilde{T}\varphi_1)(x),\ (\widetilde{T}\varphi_2)(x))$$

$$= \max_{i=1,\cdots,n}\left\{\sup_{x\in I_i} d_\infty((\widetilde{T}\varphi_1)(x),\ (\widetilde{T}\varphi_2)(x))\right\}$$

$$= \max_{i=1,\cdots,n}\left\{\sup_{x\in I_i} d_\infty(\widetilde{\mathbf{F}}_i(\widetilde{l}_i^{-1}(x),\ \varphi_1(\widetilde{l}_i^{-1}(x))),\ \widetilde{\mathbf{F}}_i(\widetilde{l}_i^{-1}(x),\ \varphi_2(\widetilde{l}_i^{-1}(x)))\right\}$$

$$= \max_{i=1,\cdots,n}\left\{\sup_{x\in I_i} d_\infty(\alpha_i\varphi_1(\widetilde{l}_i^{-1}(x))\oplus\widetilde{\mathbf{q}}_i(x),\ \alpha_i\varphi_2(\widetilde{l}_i^{-1}(x))\oplus\widetilde{\mathbf{q}}_i(x))\right\}$$

$$= \max_{i=1,\cdots,n}\left\{\sup_{x\in I_i} \alpha_i\, d_\infty(\varphi_1(\widetilde{l}_i^{-1}(x)),\ \varphi_2(\widetilde{l}_i^{-1}(x)))\right\} \le \alpha\, D(\varphi_1,\ \varphi_2),$$

where $\alpha = \max_{i=1,\cdots,n}\{\alpha_i\}(<1)$, i.e., $\widetilde{T}$ is a contractive mapping with contractive coefficient $\alpha$.

Therefore, by Banach fixed point principle, $\widetilde{T}$ has a fixed point, $\widetilde{\mathbf{f}}$ in $C^*(I,\ \mathbf{R_F})$:

$$\widetilde{\mathbf{f}}(x) = \widetilde{T}\widetilde{\mathbf{f}}(x) = \widetilde{\mathbf{F}}_i(\widetilde{l}_i^{-1}(x),\ \widetilde{\mathbf{f}}(\widetilde{l}_i^{-1}(x))) = \alpha_i\widetilde{\mathbf{f}}(\widetilde{l}_i^{-1}(x))\oplus\widetilde{\mathbf{q}}_i(x). \tag{7}$$

We show that $Gr\,\widetilde{\mathbf{f}} = A$, where $Gr\,\widetilde{\mathbf{f}} = \{(x,\ \widetilde{\mathbf{f}}(x))\in I\times\mathbf{R_F} :\ x\in I\}$. In fact, we have

$$A_i = Gr\,\widetilde{\mathbf{f}}\big|_{I_i} = \{(x,\ \widetilde{\mathbf{f}}(x))\in I\times\mathbf{R_F} :\ x\in I_i\}$$

$$= \{(x,\ \widetilde{\mathbf{F}}_i(\widetilde{l}_i^{-1}(x),\ \widetilde{\mathbf{f}}(\widetilde{l}_i^{-1}(x))))\in I\times\mathbf{R_F} :\ x\in I_i\}$$

$$= \{(\widetilde{l}_i(\widetilde{x}),\ \widetilde{\mathbf{F}}_i(\widetilde{x},\ \widetilde{\mathbf{f}}(\widetilde{x})))\in I\times\mathbf{R_F} :\ \widetilde{x}\in I_j,\ j\in\Lambda(i)\}$$

$$= \{\widetilde{\mathbf{w}}_i(\widetilde{x},\widetilde{\mathbf{f}}(\widetilde{x})) :\ \widetilde{x}\in I_j,\ j\in\Lambda(i)\}$$

$$= \bigcup_{j\in\Lambda(i)} \widetilde{\mathbf{w}}_i\left(Gr\,\widetilde{\mathbf{f}}\big|_{I_j}\right) = \bigcup_{j\in\Lambda(i)}\widetilde{\mathbf{w}}_i(A_j).$$

Hence, $Gr\,\widetilde{\mathbf{f}}$ is the attractor of RIFS (6) and since the attractor of RIFS is unique, $Gr\,\widetilde{\mathbf{f}} = A$.

□

As $\widetilde{\mathbf{f}}$, the fixed point is a continuous function interpolating the given data set $\mathbf{P}$, it is a fuzzy valued RFIF. Since $\widetilde{\mathbf{f}}$ is a fixed point of a contractive mapping $\widetilde{T}$, for any $\mathbf{f}_0\in C^*(I,\ \mathbf{R_F})$, the iterated sequence $\{\widetilde{T}^n(\mathbf{f}_0)\}$ converges to $\widetilde{\mathbf{f}}$, where $\widetilde{T}^n(\mathbf{f}_0) = \widetilde{T}(\widetilde{T}^{n-1}(\mathbf{f}_0))$ and $\widetilde{T}^0(\mathbf{f}_0) = \mathbf{f}_0$.

**Example 2.** $\mathbf{P}=\left\{\left(\frac{i}{4},\ u_i\right)\in\mathbf{R}\times\mathbf{R_F}\ \middle|\ i=0,\ 1,\cdots,\ 4\right\}$ is given, where $u_i$ s are the following triangular fuzzy numbers:

$u_0(y)=(2,\ 2,\ 2)$, $u_1(y)=(3,\ 1,\ 1)$, $u_2(y)=(5,\ 3,\ 3)$, $u_3(y)=(4,\ 2,\ 2)$, $u_4(y)=(5,\ 1,\ 1)$.

We choose the vertical scaling factors $\alpha_i$ satisfying the conditions (a1)-(a3) as follows:

$$\alpha_1=0.3,\ \alpha_2=0.33,\ \alpha_3=0.65,\ \alpha_4=0.5.$$

Let's $\widetilde{I}_1=[x_0,\ x_2],\ \widetilde{I}_2=[x_1,\ x_4],\ \widetilde{I}_3=[x_0,\ x_2],\ \widetilde{I}_4=[x_1,\ x_3]$. Then, the row stochastic matrix is as follows:

$$M=\begin{pmatrix}1/2 & 0 & 1/2 & 0\\ 1/4 & 1/4 & 1/4 & 1/4\\ 0 & 1/2 & 0 & 1/2\\ 0 & 1 & 0 & 0\end{pmatrix}.$$

The mappings $\widetilde{\mathbf{q}}_i : I_i\to\mathbf{R_F}$, $i=1,\ 2,\ 3,\ 4$ are defined as (4) in Example 1 and by (5), these are Lipschitz mappings with the following Lipschitz coefficients:

$$L_{\widetilde{\mathbf{q}}_1}=6,\ L_{\widetilde{\mathbf{q}}_2}=16.04,\ L_{\widetilde{\mathbf{q}}_3}=18.6,\ L_{\widetilde{\mathbf{q}}_4}=8.$$

The graph of a fuzzy valued RFIF $\widetilde{\mathbf{f}}$, whose graph is the attractor of RIFS (6), is as the following:

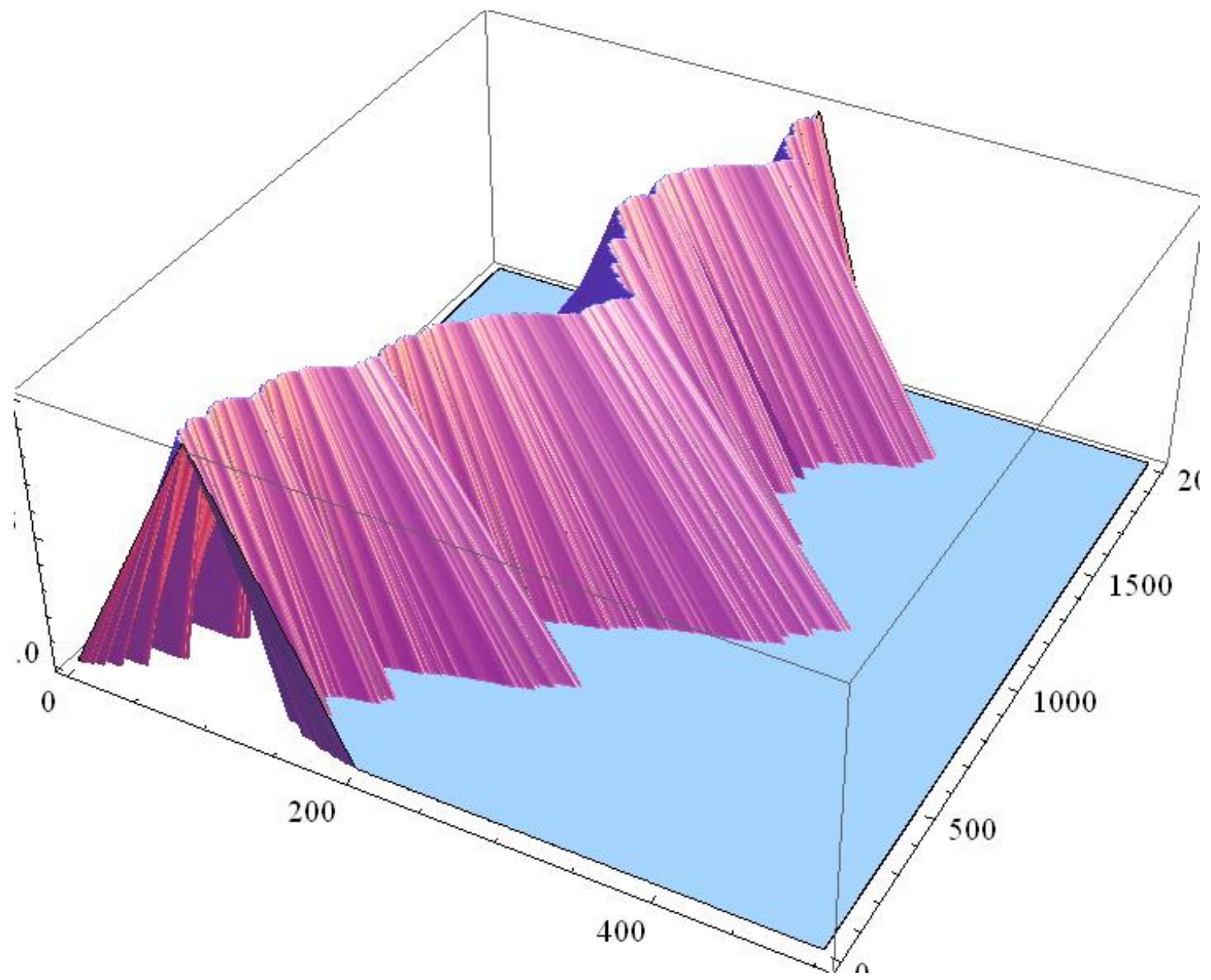

Fig 1. The graph of fuzzy valued RFIF $\widetilde{\mathbf{f}}$

The following figures show the 1-level, 0.75-level and 0.5-level set of $\widetilde{\mathbf{f}}$.

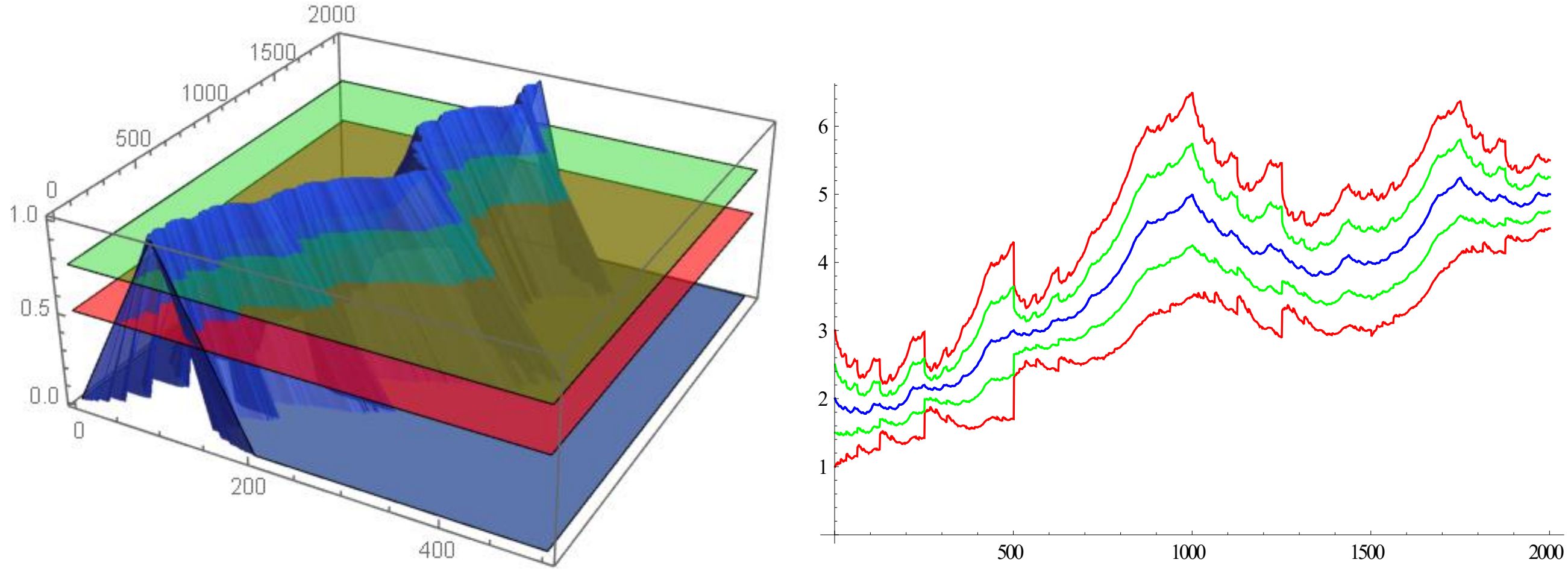

Fig 2. The 1-level set and 0.75-level set, 0.5-level set of $\tilde{\mathbf{f}}$

(Red- $(\tilde{\mathbf{f}}(x))^-(0.5)$, $(\tilde{\mathbf{f}}(x))^+(0.5)$, Green- $(\tilde{\mathbf{f}}(x))^-(0.75)$, $(\tilde{\mathbf{f}}(x))^+(0.75)$, Blue- $(\tilde{\mathbf{f}}(x))^+(1)=(\tilde{\mathbf{f}}(x))^-(1)$ )

## 3. Properties of a fuzzy valued recurrent fractal interpolaion function

In this section, we consider Hölder continuity and stability of the fuzzy valued RFIF for perturbation of the data set or the vertical scaling factors.

Let's denote the following notations:

$$L_{\tilde{\mathbf{q}}}=\max_{i=1,\cdots,n}\{L_{\tilde{\mathbf{q}}_i}\},\quad c_{\max}=\max_{i=1,\cdots,n}\{c_{\tilde{l}_i}\},\quad c_{\min}=\min_{i=1,\cdots,n}\{c_{\tilde{l}_i}\},\quad \alpha=\max_{i=1,\cdots,n}\{\alpha_i\},$$

$$|\tilde{I}|_{\max}=\max_{i=1,\cdots,n}\{|\tilde{I}_{\sigma(i)}|\},\ |\tilde{I}|_{\min}=\min_{i=1,\cdots,n}\{|\tilde{I}_{\sigma(i)}|\}.$$

### 1) Hölder continuity of a fuzzy valued RFIF

***Theorem 3. The fuzzy valued RFIF*** $\tilde{\mathbf{f}}$ ***is Hölder continuous.***

Proof. Let's $\tilde{I}_{r_1r_2\cdots r_m}=\tilde{l}_{r_m}\circ\tilde{l}_{r_{m-1}}\circ\cdots\circ\tilde{l}_{r_1}(\tilde{I}_{\sigma(i)}),\ r_i\in\{1,2,\cdots,n\}$, $L^0(X)=\{x_i \mid i=0,1,\cdots,n\}$ and $L^k(X)=\bigcup_{i=1}^{n}\{\tilde{l}_i(x)\mid x\in L^{k-1}(X)\cap\tilde{I}_{\sigma(i)}\}$ for $k=1,2,\cdots$. For any $x,\ y\in I\ (x<y)$, there exists a natural number $s$ such that

$$(x,\ y)\cap L^{s-1}(X)=\phi,\ (x,\ y)\cap L^{s}(X)\neq\phi .$$

Let's $x'\in(x,\ y)\cap L^{s}(X)$. Then there is a natural number such that

$$I_{r_1r_2\cdots r_m}\subset[x,\ x']\subset I_{r_2\cdots r_m}.$$

We have

$$d_\infty(\widetilde{\mathbf{f}}(x),\ \widetilde{\mathbf{f}}(x')) = d_\infty(\alpha_{r_m}\widetilde{\mathbf{f}}(\widetilde{l}_{r_m}^{-1}(x)) \oplus \widetilde{\mathbf{q}}_{r_m}(x),\ \alpha_{r_m}\widetilde{\mathbf{f}}(\widetilde{l}_{r_m}^{-1}(x')) \oplus \widetilde{\mathbf{q}}_{r_m}(x'))$$

$$\le \alpha_{r_m} d_\infty(\widetilde{\mathbf{f}}(\widetilde{l}_{r_m}^{-1}(x)),\ \widetilde{\mathbf{f}}(\widetilde{l}_{r_m}^{-1}(x'))) + d_\infty(\widetilde{\mathbf{q}}_{r_m}(x),\ \widetilde{\mathbf{q}}_{r_m}(x'))$$

$$\le \alpha\, d_\infty(\widetilde{\mathbf{f}}(\widetilde{l}_{r_m}^{-1}(x)),\ \widetilde{\mathbf{f}}(\widetilde{l}_{r_m}^{-1}(x'))) + L_{\widetilde{\mathbf{q}}_{r_m}}\ |x-x'|$$

$$\le \alpha[\alpha\, d_\infty(\widetilde{\mathbf{f}}(\widetilde{l}_{r_{m-1}}^{-1}\circ\widetilde{l}_{r_m}^{-1}(x)),\ \widetilde{\mathbf{f}}(\widetilde{l}_{r_{m-1}}^{-1}\circ\widetilde{l}_{r_m}^{-1}(x'))) + L_{\widetilde{\mathbf{q}}} c_{r_m}^{-1}\,|x-x'|] + L_{\widetilde{\mathbf{q}}}\,|x-x'|$$

$$\le \alpha^2\, d_\infty(\widetilde{\mathbf{f}}(\widetilde{l}_{r_{m-1}}^{-1}\circ\widetilde{l}_{r_m}^{-1}(x)),\ \widetilde{\mathbf{f}}(\widetilde{l}_{r_{m-1}}^{-1}\circ\widetilde{l}_{r_m}^{-1}(x'))) + L_{\widetilde{\mathbf{q}}}\left(1+\frac{\alpha}{c_{\min}}\right)\cdot|x-x'|$$

$$\le \cdots \le$$

$$\le \alpha^{m-1}\, d_\infty(\widetilde{\mathbf{f}}(\widetilde{l}_{r_2}^{-1}\circ\cdots\circ\widetilde{l}_{r_m}^{-1}(x)),\ \widetilde{\mathbf{f}}(\widetilde{l}_{r_2}^{-1}\circ\cdots\circ\widetilde{l}_{r_m}^{-1}(x'))) + L_{\widetilde{\mathbf{q}}}\left(1+\frac{\alpha}{c_{\min}}+\cdots\left(\frac{\alpha}{c_{\min}}\right)^{m-2}\right)\cdot|x-x'|.$$

For $\mathbf{0}$, fuzzy zero number, we get

$$\begin{aligned} d_\infty(\widetilde{\mathbf{f}}(x),\ \mathbf{0}) &= d_\infty(\alpha_i\widetilde{\mathbf{f}}(\widetilde{l}_i^{-1}(x)) \oplus \widetilde{\mathbf{q}}_i(x),\ \mathbf{0}) \\ &\le \alpha_i d_\infty(\widetilde{\mathbf{f}}(\widetilde{l}_i^{-1}(x)),\ \mathbf{0}) + d_\infty(\widetilde{\mathbf{q}}_i(x),\ \mathbf{0}) \\ &\le \alpha\, D(\widetilde{\mathbf{f}},\ \mathbf{0}) + L_{\widetilde{\mathbf{q}}_i}\cdot|x| \end{aligned}$$

and hence, we have $D(\widetilde{\mathbf{f}},\ \mathbf{0}) \le \alpha\, D(\widetilde{\mathbf{f}},\ \mathbf{0}) + L_{\widetilde{\mathbf{q}}}\cdot\max\{|x_0|,\ |x_n|\}$. As a result,

$$D(\widetilde{\mathbf{f}},\ \mathbf{0}) \le \frac{L_{\widetilde{\mathbf{q}}}\cdot\max\{|x_0|,\ |x_n|\}}{1-\alpha}. \tag{7}$$

Let's $M = \dfrac{2L_{\widetilde{\mathbf{q}}}\cdot\max\{|x_0|,\ |x_n|\}}{1-\alpha}$. Then we get

$$d_\infty(\widetilde{\mathbf{f}}(\widetilde{l}_{r_2}^{-1}\circ\cdots\circ\widetilde{l}_{r_m}^{-1}(x)),\ \widetilde{\mathbf{f}}(\widetilde{l}_{r_2}^{-1}\circ\cdots\circ\widetilde{l}_{r_m}^{-1}(x'))) \le 2D(\widetilde{\mathbf{f}},\ \mathbf{0}) \le M$$

and

$$d_\infty(\widetilde{\mathbf{f}}(x),\ \widetilde{\mathbf{f}}(x')) \le M\alpha^{m-1} + L_{\widetilde{\mathbf{q}}}\left(1+\frac{\alpha}{c_{\min}}+\cdots\left(\frac{\alpha}{c_{\min}}\right)^{m-2}\right)\cdot|x-x'|.$$

We denote $\delta = \alpha / c_{\min}$. As $I_{r_1 r_2 \ldots r_m} \subset [x,\ x'] \subset I_{r_2\cdots r_m}$, we obtain

$$c_{\min}^m\,|\widetilde{I}_{\sigma(i)}| \le |x-x'| \le c_{\max}^{m-1}\,|\widetilde{I}_{\sigma(i)}|$$

and

$$d_\infty(\widetilde{\mathbf{f}}(x),\ \widetilde{\mathbf{f}}(x')) \le M\delta^{m-1}c_{\min}^{m-1} + L_{\widetilde{\mathbf{q}}}(1+\delta+\cdots\delta^{m-2})\cdot|x-x'|$$
$$\le M\frac{|x-x'|}{c_{\min}|\widetilde{I}|_{\min}}\delta^{m-1} + L_{\widetilde{\mathbf{q}}}(1+\delta+\cdots\delta^{m-2})\cdot|x-x'| \le N\sum_{j=0}^{m-1}\delta^j\cdot|x-x'|,$$

where $N=\max\left\{\frac{M}{c_{\min}|\widetilde{I}|_{\min}},\ L_{\widetilde{\mathbf{q}}}\right\}$.

In [18], they showed that $N\sum_{j=0}^{m-1}\delta^j\cdot|x-x'| \le Q|x-x'|^\tau$, where

(1) in the case of $\delta<1$, i.e. $\alpha<c_{\min}$, $\tau=1$, $Q=\frac{N}{1-\delta}$,

(2) in the case of $\delta=1$, i.e. $\alpha=c_{\min}$, $0<\tau<1$, $Q=N\left(1-\frac{1}{(1-\tau)\mathrm{e}\ln c_{\max}}\right)\max\{1,\ |\widetilde{I}|_{\max}\}$,

(3) in the case of $\delta>1$, i.e. $\alpha>c_{\min}$, $\tau=\frac{\ln\delta}{\ln c_{\max}}+1$, $Q=\frac{N\delta}{\delta-1}\max\{1,\ |\widetilde{I}|_{\max}\}$.

Hence, we get

$$d_\infty(\widetilde{\mathbf{f}}(x),\ \widetilde{\mathbf{f}}(x')) \le Q|x-x'|^\tau$$

and similarly, we have

$$d_\infty(\widetilde{\mathbf{f}}(y),\ \widetilde{\mathbf{f}}(x')) \le Q|y-x'|^\tau.$$

As a result, we obtain

$$d_\infty(\widetilde{\mathbf{f}}(x),\ \widetilde{\mathbf{f}}(y)) \le Q|x-x'|^\tau + Q|y-x'|^\tau \le 2Q|x-y|^\tau.\ \square$$

For example, for a fuzzy valued RFIF $\widetilde{\mathbf{f}}$ constructed in Example 2, since $\alpha=0.65,\ c_{\min}=1/3$, $\delta=\frac{\alpha}{c_{\min}}=\frac{0.65}{1/3}>1$ and by (3), the Hölder exponent of $\widetilde{\mathbf{f}}$ is as follows:

$$\tau=\frac{\ln\delta}{\ln c_{\max}}+1=\ln\left(\frac{0.65}{1/3}\right)\Big/\ln 0.5+1\approx 0.0365.$$

**2) Stability of a fuzzy valued RFIF for perturbation of the data set**

We show that any small perturbation in the data set results in only small perturbation of fuzzy valued RFIF. To do it, we study stability of a fuzzy valued RFIF for perturbation of each coordinate of points in the data set and then, for perturbation of all coordinates.

Firstly, let's $\widetilde{\mathbf{f}}^{x*}$ be a fuzzy valued RFIF interpolating $\mathbf{P}^{x*}=\{(x_i^*,\ y_i):\ i=0,\ \cdots,\ n\}$, in which $x$-coordinates of points are perturbed as follows:

$$x_0=x_0^*<x_1^*<\cdots<x_n^*=x_n\ .$$

Then,

$$\widetilde{\mathbf{f}}^{x*}(x)=\alpha_i\widetilde{\mathbf{f}}^{x*}(\widetilde{l}_i^{*-1}(x))\oplus\widetilde{\mathbf{q}}_i^{x*}(x),\ \ x\in I_i^*,$$

where $I_i^*=[x_{i-1}^*,x_i^*]$ and $\widetilde{l}_i^*=R\circ\widetilde{l}_i\circ R^{-1}$, $\widetilde{\mathbf{q}}_i^{x*}=\widetilde{\mathbf{q}}_i\circ R^{-1}$ for the following function $R:I_i\to I_i^*$:

$$R(x)=x_{i-1}^*+\frac{x_i^*-x_{i-1}^*}{x_i-x_{i-1}}(x-x_{i-1}),\ \ x\in I_i\,.$$

We can easily show that $\widetilde{\mathbf{q}}_i^{x*}:\ I_i^*\to\mathbf{R}_{\mathbf{F}}$ is a Lipschitz mapping satisfying the following condition:

$$\alpha_i u_{s_i}\oplus\widetilde{\mathbf{q}}_i^{x*}(\widetilde{l}_i^*(x_{s_i}^*))=u_{i-1},\ \ \alpha_i u_{e_i}\oplus\widetilde{\mathbf{q}}_i^{x*}(\widetilde{l}_i^*(x_{e_i}^*))=u_i\,.$$

The function $R(x)$ satisfies the following inequality ([22]):

$$|\,R(x)-x\,|\le\max_{i=0,\cdots,n}|\,x_i-x_i^*\,|,\ \ x\in I\ .$$

***Lemma 1. The fuzzy valued RFIFs $\widetilde{\mathbf{f}}$ and $\widetilde{\mathbf{f}}^{x*}$ interpolating $\mathbf{P}$ and $\mathbf{P}^{x*}$, respectively, satisfy the following:***

$$D(\widetilde{\mathbf{f}},\ \widetilde{\mathbf{f}}^{x*})\le\frac{(1+\alpha)H_{\widetilde{\mathbf{f}}}}{1-\alpha}\max_{i=0,\cdots,n}|\,x_i-x_i^*\,|^{\tau},$$

***where $\tau$ is the Hölder exponent of $\widetilde{\mathbf{f}}$.***

Proof. For any $x\in I_i^*$, we have

$$\begin{aligned}d_\infty(\widetilde{\mathbf{f}}(x),\ \widetilde{\mathbf{f}}^{x*}(x))&\le d_\infty(\widetilde{\mathbf{f}}(x),\ \widetilde{\mathbf{f}}(R^{-1}(x)))+d_\infty(\widetilde{\mathbf{f}}(R^{-1}(x)),\ \widetilde{\mathbf{f}}^{x*}(x))\\&\le H_{\widetilde{\mathbf{f}}}\,|\,x-R^{-1}(x)\,|^{\tau}+d_\infty(\alpha_i\widetilde{\mathbf{f}}(\widetilde{l}_i^{-1}(R^{-1}(x)))\oplus\widetilde{\mathbf{q}}_i(R^{-1}(x)),\ \alpha_i\widetilde{\mathbf{f}}^{x*}(\widetilde{l}_i^{*-1}(x))\oplus\widetilde{\mathbf{q}}_i^{x*}(x))\\&\le H_{\widetilde{\mathbf{f}}}\max_{i=0,\cdots,n}|\,x_i-x_i^*\,|^{\tau}+\alpha_i\,d_\infty(\widetilde{\mathbf{f}}(\widetilde{l}_i^{-1}(R^{-1}(x))),\ \widetilde{\mathbf{f}}^{x*}(\widetilde{l}_i^{*-1}(x)))\end{aligned}$$

$$\le H_{\tilde{\mathbf{f}}} \max_{i=0,\cdots,n} |x_i - x_i^*|^\tau + \alpha_i\, d_\infty(\tilde{\mathbf{f}}(\tilde{l}_i^{-1} \circ R^{-1}(x)),\ \tilde{\mathbf{f}}(R \circ \tilde{l}_i^{-1} \circ R^{-1}(x)))$$

$$+ \alpha_i\, d_\infty(\tilde{\mathbf{f}}(R \circ \tilde{l}_i^{-1} \circ R^{-1}(x)),\ \tilde{\mathbf{f}}^{x*}(R \circ \tilde{l}_i^{-1} \circ R^{-1}(x)))$$

$$\le (1+\alpha_i) H_{\tilde{\mathbf{f}}} \max_{i=0,\cdots,n} |x_i - x_i^*|^\tau + \alpha_i\, D(\tilde{\mathbf{f}},\ \tilde{\mathbf{f}}^{x*})$$

and we get

$$D(\tilde{\mathbf{f}},\ \tilde{\mathbf{f}}^{x*}) \le (1+\alpha) H_{\tilde{\mathbf{f}}} \max_{i=0,\cdots,n} |x_i - x_i^*|^\tau + \alpha\, D(\tilde{\mathbf{f}},\ \tilde{\mathbf{f}}^{x*}),$$

$$D(\tilde{\mathbf{f}},\ \tilde{\mathbf{f}}^{x*}) \le \frac{(1+\alpha) H_{\tilde{\mathbf{f}}}}{1-\alpha} \max_{i=0,\cdots,n} |x_i - x_i^*|^\tau. \ \square$$

Secondly, let's $\tilde{\mathbf{f}}^{u*}$ be a fuzzy valued RFIF interpolating $\mathbf{P}^{u*} = \{(x_i,\ u_i^*) : i = 0,\ \cdots,\ n\}$, in which $u$-coordinates of points are perturbed, then,

$$\tilde{\mathbf{f}}^{u*}(x) = \alpha_i \tilde{\mathbf{f}}^{u*}(\tilde{l}_i^{-1}(x)) \oplus \tilde{\mathbf{q}}_i^{u*}(x),\ x \in I_i,$$

where $\tilde{\mathbf{q}}_i^{u*},\ i = 1,\ \cdots,\ n$ are Lipschitz mappings satisfying

$$\alpha_i u_{s_i}^* \oplus \tilde{\mathbf{q}}_i^{u*}(\tilde{l}_i(x_{s_i})) = u_{i-1}^*,\ \alpha_i u_{e_i}^* \oplus \tilde{\mathbf{q}}_i^{u*}(\tilde{l}_i(x_{e_i})) = u_i^*$$

and

$$D(\tilde{\mathbf{q}}_i,\ \tilde{\mathbf{q}}_i^{u*}) \le \mu \max_{i=0,\cdots,n} \{d_\infty(u_i,\ u_i^*)\}.$$

***Lemma 2. The fuzzy valued RFIFs $\tilde{\mathbf{f}}$ and $\tilde{\mathbf{f}}^{u*}$ interpolating $\mathbf{P}$ and $\mathbf{P}^{u*}$, respectively, satisfy the following:***

$$D(\tilde{\mathbf{f}},\ \tilde{\mathbf{f}}^{u*}) \le \frac{\mu}{1-\alpha} \max_{i=0,\cdots,n} \{d_\infty(u_i,\ u_i^*)\}.$$

Proof. For any $x \in I_i$, we have

$$d_\infty(\tilde{\mathbf{f}}(x),\ \tilde{\mathbf{f}}^{u*}(x)) \le d_\infty(\alpha_i \tilde{\mathbf{f}}(\tilde{l}_i^{-1}(x)) \oplus \tilde{\mathbf{q}}_i(x),\ \alpha_i \tilde{\mathbf{f}}^{u*}(\tilde{l}_i^{-1}(x)) \oplus \tilde{\mathbf{q}}_i^{u*}(x))$$

$$\le \alpha_i\, d_\infty(\tilde{\mathbf{f}}(\tilde{l}_i^{-1}(x)),\ \tilde{\mathbf{f}}^{u*}(\tilde{l}_i^{-1}(x))) + d_\infty(\tilde{\mathbf{q}}_i(x),\ \tilde{\mathbf{q}}_i^{u*}(x))$$

$$\le \alpha_i\, D(\tilde{\mathbf{f}},\ \tilde{\mathbf{f}}^{u*}) + \mu \max_{i=0,\cdots,n} \{d_\infty(u_i,\ u_i^*)\}$$

and we have

$$D(\tilde{\mathbf{f}},\ \tilde{\mathbf{f}}^{u*}) \le \alpha\, D(\tilde{\mathbf{f}},\ \tilde{\mathbf{f}}^{u*}) + \mu \max_{i=0,\cdots,n} \{d_\infty(u_i,\ u_i^*)\},$$

$$D(\widetilde{\mathbf{f}},\ \widetilde{\mathbf{f}}^{u*}) \le \frac{\mu}{1-\alpha} \max_{i=0,\cdots,n} \{d_\infty(u_i,\ u_i^*)\}. \ \square$$

Finally, let's $\widetilde{\mathbf{f}}^*$ be a fuzzy valued RFIF interpolating the perturbed data set $\mathbf{P}^* = \{(x_i^*,\ u_i^*) :\ i = 0,\ \cdots,\ n\}$:

$$\widetilde{\mathbf{f}}^*(x) = \alpha_i \widetilde{\mathbf{f}}^*(\widetilde{l}_i^{*-1}(x)) \oplus \widetilde{\mathbf{q}}_i^*(x),\ \ x \in I_i^*,$$

where $\widetilde{\mathbf{q}}_i^*$s are Lipschitz mappings satisfying

$$\alpha_i u_{s_i}^* \oplus \widetilde{\mathbf{q}}_i^{x*}(\widetilde{l}_i^*(x_{s_i}^*)) = u_{i-1}^*,\ \ \alpha_i u_{e_i}^* \oplus \widetilde{\mathbf{q}}_i^{x*}(\widetilde{l}_i^*(x_{e_i}^*)) = u_i^*$$

and

$$D(\widetilde{\mathbf{q}}_i \circ R^{-1},\ \widetilde{\mathbf{q}}_i^*) \le \mu \max_{i=0,\cdots,n} \{d_\infty(u_i,\ u_i^*)\}.$$

From Lemma 1 and 2, we have the following result:

***Theorem 4. The fuzzy valued RFIFs*** $\widetilde{\mathbf{f}}$ ***and*** $\widetilde{\mathbf{f}}^*$ ***interpolating*** $\mathbf{P}$ ***and*** $\mathbf{P}^*$***, respectively, satisfy the following:***

$$D(\widetilde{\mathbf{f}},\ \widetilde{\mathbf{f}}^*) \le \frac{(1+\alpha) H_{\widetilde{\mathbf{f}}} \max_{i=0,\cdots,n} |x_i - x_i^*|^\tau + \mu \max_{i=0,\cdots,n} \{d_\infty(u_i,\ u_i^*)\}}{1-\alpha}.$$

**Example 3.** Let's $\widetilde{\mathbf{f}}$ and $\widetilde{\mathbf{f}}^*$ be fuzzy valued RFIF with the vertical scaling factors $\alpha_i$ satisfying the conditions (a1)-(a3) in Example 1 and $\widetilde{\mathbf{q}}_i$, $\widetilde{\mathbf{q}}_i^*$ defined in (4) with $\mathbf{P}$ and $\mathbf{P}^*$, respectively.

Since for any $x \in I_i^*$,

$$d_\infty(\widetilde{\mathbf{q}}_i \circ R^{-1}(x),\ \widetilde{\mathbf{q}}_i^*(x)) =$$

$$= d_\infty\left(\left(\frac{R^{-1}(x) - x_{i-1}}{x_i - x_{i-1}} u_i \oplus \frac{R^{-1}(x) - x_i}{x_{i-1} - x_i} u_{i-1}\right) \vee_H \left[\alpha_i\left(\frac{\widetilde{l}_i^{-1}(R^{-1}(x)) - x_{s_i}}{x_{e_i} - x_{s_i}} u_{e_i} \oplus \frac{\widetilde{l}_i^{-1}(R^{-1}(x)) - x_{e_i}}{x_{s_i} - x_{e_i}} u_{s_i}\right)\right],\right.$$

$$\left.\left(\frac{x - x_{i-1}^*}{x_i^* - x_{i-1}^*} u_i^* \oplus \frac{x - x_i^*}{x_{i-1}^* - x_i^*} u_{i-1}^*\right) \vee_H \left[\alpha_i\left(\frac{\widetilde{l}_i^{*-1}(x) - x_{s_i}^*}{x_{e_i}^* - x_{s_i}^*} u_{e_i}^* \oplus \frac{\widetilde{l}_i^{*-1}(x) - x_{e_i}^*}{x_{s_i}^* - x_{e_i}^*} u_{s_i}^*\right)\right]\right)$$

$$\le d_\infty\left(\frac{x-x_{i-1}^*}{x_i^*-x_{i-1}^*}u_i \oplus \frac{x-x_i^*}{x_{i-1}^*-x_i^*}u_{i-1},\ \frac{x-x_{i-1}^*}{x_i^*-x_{i-1}^*}u_i^* \oplus \frac{x-x_i^*}{x_{i-1}^*-x_i^*}u_{i-1}^*\right)+$$

$$+d_\infty\left(\alpha_i\left(\frac{x-x_{i-1}^*}{x_i^*-x_{i-1}^*}u_{e_i} \oplus \frac{x-x_i^*}{x_{i-1}^*-x_i^*}u_{s_i}\right),\ \alpha_i\left(\frac{x-x_{i-1}^*}{x_i^*-x_{i-1}^*}u_{e_i}^* \oplus \frac{x-x_i^*}{x_{i-1}^*-x_i^*}u_{s_i}^*\right)\right)$$

$$\le \frac{x-x_{i-1}^*}{x_i^*-x_{i-1}^*}[d_\infty(u_i,\ u_i^*)+\alpha_i d_\infty(u_{e_i},\ u_{e_i}^*)]+\frac{x_i^*-x}{x_i^*-x_{i-1}^*}[d_\infty(u_{i-1},\ u_{i-1}^*)+\alpha_i d_\infty(u_{s_i},\ u_{s_i}^*)]$$

$$\le (1+\alpha)\max_{i=0,\cdots,n}\{d_\infty(u_i,\ u_i^*)\},$$

by Theorem 4, $\widetilde{\mathbf{f}}$ and $\widetilde{\mathbf{f}}^*$ satisfy the following inequality:

$$D(\widetilde{\mathbf{f}},\ \widetilde{\mathbf{f}}^*)\le \frac{1+\alpha}{1-\alpha}(H_{\widetilde{\mathbf{f}}}\max_{i=0,\cdots,n}|x_i-x_i^*|^\tau+\max_{i=0,\cdots,n}\{d_\infty(u_i,\ u_i^*)\}).$$

**3) Stability of a fuzzy valued RFIF for perturbation in the vertical scaling factors**

We show that any small perturbation in the vertical scaling factors leads to a small perturbation in fuzzy valued RFIF.

Let's $\widetilde{\mathbf{f}}$ and $\widetilde{\mathbf{f}}^*$ be the fuzzy valued RFIFs interpolating $\mathbf{P}$ with the vertical scaling factors $\alpha_i$ and $\alpha_i^*$, respectively, then

$$\widetilde{\mathbf{f}}(x)=\alpha_i\widetilde{\mathbf{f}}(\widetilde{l}_i^{-1}(x))\oplus\widetilde{\mathbf{q}}_i(x),\quad \widetilde{\mathbf{f}}^*(x)=\alpha_i^*\widetilde{\mathbf{f}}^*(\widetilde{l}_i^{-1}(x))\oplus\widetilde{\mathbf{q}}_i^*(x),\ x\in I_i,$$

where $\widetilde{\mathbf{q}}_i$ and $\widetilde{\mathbf{q}}_i^*$ are Lipschitz mappings satisfying

$$D(\widetilde{\mathbf{q}}_i,\ \widetilde{\mathbf{q}}_i^*)\le \mu\max_{i=1,\cdots,n}\{|\alpha_i-\alpha_i^*|\}.$$

Let's $\alpha^*=\max_{i=1,\cdots,n}\{\alpha_i^*\}$.

***Theorem 5. The fuzzy valued RFIFs*** $\widetilde{\mathbf{f}}$ ***and*** $\widetilde{\mathbf{f}}^*$ ***with the vertical scaling factors*** $\alpha_i$ ***and*** $\alpha_i^*$***, respectively, satisfy the following:***

$$D(\widetilde{\mathbf{f}},\ \widetilde{\mathbf{f}}^*)\le\left[\frac{L_{\widetilde{\mathbf{q}}}\cdot\max\{|x_0|,\ |x_n|\}}{(1-\alpha)(1-\alpha^*)}+\frac{\mu}{1-\alpha^*}\right]\max_{i=1,\cdots,n}\{|\alpha_i-\alpha_i^*|\}.$$

Proof. For any $x\in I_i$,

$$d_\infty(\widetilde{\mathbf{f}}(x),\ \widetilde{\mathbf{f}}^*(x))\le d_\infty(\alpha_i\widetilde{\mathbf{f}}(\widetilde{l}_i^{-1}(x))\oplus\widetilde{\mathbf{q}}_i(x),\ \alpha_i^*\widetilde{\mathbf{f}}^*(\widetilde{l}_i^{-1}(x))\oplus\widetilde{\mathbf{q}}_i^*(x))$$

$$\le d_\infty(\alpha_i \widetilde{\mathbf{f}}(\widetilde{l}_i^{-1}(x)),\ \alpha_i^* \widetilde{\mathbf{f}}^*(\widetilde{l}_i^{-1}(x))) + d_\infty(\widetilde{\mathbf{q}}_i(x),\ \widetilde{\mathbf{q}}_i^*(x))$$

$$\le d_\infty(\alpha_i \widetilde{\mathbf{f}}(\widetilde{l}_i^{-1}(x)),\ \alpha_i^* \widetilde{\mathbf{f}}(\widetilde{l}_i^{-1}(x))) + d_\infty(\alpha_i^* \widetilde{\mathbf{f}}(\widetilde{l}_i^{-1}(x)),\ \alpha_i^* \widetilde{\mathbf{f}}^*(\widetilde{l}_i^{-1}(x))) + D(\widetilde{\mathbf{q}}_i,\ \widetilde{\mathbf{q}}_i^*)$$

$$\le |\alpha_i - \alpha_i^*| \cdot d_\infty(\widetilde{\mathbf{f}}(\widetilde{l}_i^{-1}(x)),\ \mathbf{0}) + \alpha_i^* d_\infty(\widetilde{\mathbf{f}}(\widetilde{l}_i^{-1}(x)),\ \widetilde{\mathbf{f}}^*(\widetilde{l}_i^{-1}(x))) + \mu \max_{i=1,\cdots,n}\{|\alpha_i - \alpha_i^*|\}$$

$$\le [D(\widetilde{\mathbf{f}},\ \mathbf{0}) + \mu] \max_{i=1,\cdots,n}\{|\alpha_i - \alpha_i^*|\} + \alpha_i^* D(\widetilde{\mathbf{f}},\ \widetilde{\mathbf{f}}^*)$$

$$\le \left[\frac{L_{\widetilde{\mathbf{q}}} \cdot \max\{|x_0|,\ |x_n|\}}{1-\alpha} + \mu\right] \max_{i=1,\cdots,n}\{|\alpha_i - \alpha_i^*|\} + \alpha^* D(\widetilde{\mathbf{f}},\ \widetilde{\mathbf{f}}^*).$$

Hence, we get

$$D(\widetilde{\mathbf{f}},\ \widetilde{\mathbf{f}}^*) \le \left[\frac{L_{\widetilde{\mathbf{q}}} \cdot \max\{|x_0|,\ |x_n|\}}{1-\alpha} + \mu\right] \max_{i=1,\cdots,n}\{|\alpha_i - \alpha_i^*|\} + \alpha^* D(\widetilde{\mathbf{f}},\ \widetilde{\mathbf{f}}^*)$$

and the proof has been done. □

**Example 4.** Let's $\widetilde{\mathbf{f}}$ and $\widetilde{\mathbf{f}}^*$ be fuzzy valued RFIFs interpolating the same data set $\mathbf{P}$ with the different vertical scaling factors $\alpha_i$ and $\alpha_i^*$ and $\widetilde{\mathbf{q}}_i$, $\widetilde{\mathbf{q}}_i^*$ defined by (4) in Example 1, respectively. Now, we estimate error between them.

For any $x \in I_i$,

$$d_\infty(\widetilde{\mathbf{q}}_i(x),\ \widetilde{\mathbf{q}}_i^*(x)) = d_\infty\left(\left(\frac{x - x_{i-1}}{x_i - x_{i-1}} u_i \oplus \frac{x - x_i}{x_{i-1} - x_i} u_{i-1}\right) \vee_H \left[\alpha_i\left(\frac{\widetilde{l}_i^{-1}(x) - x_{s_i}}{x_{e_i} - x_{s_i}} u_{e_i} \oplus \frac{\widetilde{l}_i^{-1}(x) - x_{e_i}}{x_{s_i} - x_{e_i}} u_{s_i}\right)\right],\right.$$

$$\left.\left(\frac{x - x_{i-1}}{x_i - x_{i-1}} u_i \oplus \frac{x - x_i}{x_{i-1} - x_i} u_{i-1}\right) \vee_H \left[\alpha_i^*\left(\frac{\widetilde{l}_i^{-1}(x) - x_{s_i}}{x_{e_i} - x_{s_i}} u_{e_i} \oplus \frac{\widetilde{l}_i^{-1}(x) - x_{e_i}}{x_{s_i} - x_{e_i}} u_{s_i}\right)\right]\right)$$

$$\le d_\infty\left(\alpha_i\left(\frac{x - x_{i-1}}{x_i - x_{i-1}} u_{e_i} \oplus \frac{x - x_i}{x_{i-1} - x_i} u_{s_i}\right),\ \alpha_i^*\left(\frac{x - x_{i-1}}{x_i - x_{i-1}} u_{e_i} \oplus \frac{x - x_i}{x_{i-1} - x_i} u_{s_i}\right)\right)$$

$$\le |\alpha_i - \alpha_i^*|\, d_\infty\left(\frac{x - x_{i-1}}{x_i - x_{i-1}} u_{e_i} \oplus \frac{x - x_i}{x_{i-1} - x_i} u_{s_i},\ \mathbf{0}\right)$$

$$\le |\alpha_i - \alpha_i^*| \left[\frac{x - x_{i-1}}{x_i - x_{i-1}} d_\infty(u_{e_i},\ \mathbf{0}) + \frac{x - x_i}{x_i - x_{i-1}} d_\infty(u_{s_i},\ \mathbf{0})\right]$$

$$\le |\alpha_i - \alpha_i^*| \left[ \frac{x - x_{i-1}}{x_i - x_{i-1}} d_\infty(u_{e_i}, \mathbf{0}) + \frac{x - x_i}{x_i - x_{i-1}} d_\infty(u_{s_i}, \mathbf{0}) \right]$$

$$\le \max_{i=0,\cdots,n} \{d_\infty(u_i, \mathbf{0})\} \cdot \max_{i=1,\cdots,n} \{|\alpha_i - \alpha_i^*|\}.$$

By Theorem 5, the fuzzy valued RFIFs $\tilde{\mathbf{f}}$ and $\tilde{\mathbf{f}}^*$ satisfy the following:

$$D(\tilde{\mathbf{f}}, \tilde{\mathbf{f}}^*) \le \left[ \frac{L_{\tilde{\mathbf{q}}} \cdot \max\{|x_0|, |x_n|\}}{(1-\alpha)(1-\alpha^*)} + \frac{\max_{i=0,\cdots,n} \{d_\infty(u_i, \mathbf{0})\}}{1-\alpha^*} \right] \max_{i=1,\cdots,n} \{|\alpha_i - \alpha_i^*|\}.$$

## 4. Conclusion

In this paper, firstly, we proposed a construction of fuzzy valued RFIFs using RIFSs as interpolation functions that are very useful to model objects with both local self-similarity and uncertainty. Next, we studied Hölder continuity and stability of the interpolation functions due to perturbations of the data set or the vertical scaling factors.

In the future, we are going to study fuzzy calculus and fractional calculus of fuzzy valued RFIFs using their properties.